\documentclass[11pt,dvipsnames]{amsart}
\usepackage{amsmath,amssymb,amsfonts}
\usepackage[mathscr]{eucal}
\usepackage{mathtools}

\usepackage{tikz}%
\usetikzlibrary{matrix,arrows,decorations.pathmorphing,cd,patterns,calc,backgrounds,patterns}
\tikzset{dummy/.style= {circle,fill,draw,inner sep=0pt,minimum size=1.2mm}}
\tikzset{vertex/.style={fill, circle, minimum size=.1cm, inner sep=0pt}}

\usepackage{multirow}

\usepackage[colorlinks, citecolor=PineGreen, linkcolor=PineGreen]{hyperref}
\hypersetup{linktoc=all} 

\usepackage[textwidth=25mm, textsize=tiny]{todonotes}
\usepackage[final]{microtype}

\usepackage{enumerate}

\numberwithin{equation}{section} 
\numberwithin{figure}{section}
\usepackage[nameinlink,capitalise,noabbrev]{cleveref}

\newcommand{\newrefformat}[2]{}

\crefname{lemma}{el Lema}{los Lemas}
\crefname{thm}{el Teorema}{los Teorema}
\crefname{defn}{la Definición}{las Definiciones}
\crefname{propn}{la Proposición}{las Proposiciones}
\crefname{obs}{la Observación}{las Observaciones}
\crefname{cor}{el Corolario}{los Corolarios}
\crefname{equation}{la Ecuación}{las Ecuaciones}
\crefname{notation}{la Notación}{las Notaciones}
\crefname{const}{la Construcción}{las Construcciones}
\crefname{ex}{el Ejemplo}{los Ejemplos}

\theoremstyle{plain}
\newtheorem{thm}[equation]{Teorema}
\newtheorem{cor}[equation]{Corolario}
\newtheorem{propn}[equation]{Proposición}
\newtheorem{lemma}[equation]{Lema}

\theoremstyle{definition}
\newtheorem{defn}[equation]{Definición}
\newtheorem{ex}[equation]{Ejemplo}

\newtheorem{notation}[equation]{Notación}

\usepackage[style = alphabetic, doi=false]{biblatex}
\usepackage{hyperref}
\usepackage{xurl}
\hypersetup{breaklinks=true}
\renewbibmacro{in:}{}
\definecolor{winered}{rgb}{0.5,0,0}
 \usepackage{hyperref}
 \hypersetup{citecolor = winered, linkcolor = winered, menucolor = winered, colorlinks = true}

\newcommand{\mapdel}{\Delta}

\renewcommand{\phi}{\varphi}

\newcommand{\Sing}{\operatorname{Sing}}

\input xy
\xyoption{all}

\title{Teoría de homotopía usando Conjuntos Simpliciales}

\author[M. Campillo]{Mathilda Campillo} 
\address[M. Campillo]{Departamento de Matemáticas y Estadística, Universidad del Norte,
Km. 5 Vía Antigua Puerto Colombia, Barranquilla 081007, Colombia}
\email{mathildac@uninorte.edu.co}

\author[A. Osorno]{Angélica M.~Osorno}
\address[A. Osorno]{Department of Mathematics and Statistics, Reed College, Portland, OR 97202, USA}
\email{aosorno@reed.edu}

\author[M. Rivera]{Manuel Rivera}
\address[M. Rivera]{Department of Mathematics, Purdue University, 150 N University St.,
West Lafayette, IN, 47906}
\email{manuelr@purdue.edu}


\subjclass[2020]{}

\begin{document}

\maketitle
\section{Introducción}
Estas notas surgieron como parte de un mini-curso dictado por A. Osorno y M. Rivera en el primer \textit{Encuentro Colombiano de Geometría y Topología} (ECOGyT) que se llevó a cabo en la Universidad Nacional de Colombia sede Bogotá en julio del 2024. La idea es que sirvan como una guía para un primer encuentro con la teoría de homotopía y técnicas simpliciales -enfatizando en la intuición y enunciados importantes - accesible a estudiantes con conocimiento básico de topología general y así invitar a indagar más profundamente sobre el tema y sus aplicaciones en distintos campos. 

No hemos incluido demostraciones en este documento. La mejor manera de utilizar estas notas es \textit{tratando de entender y demostrar cada enunciado (lemas, proposiciones, teoremas o aseveraciones hechas en medio del texto) por cuenta propia sin consultar las referencias o la literatura}. Algunos enunciados son más difíciles que otros y lo ideal sería que el lector mismo lo descubra. De completar esta tarea, invitamos al lector interesado en aprender más a contactar a alguno de los autores.

\subsection*{Agradecimientos} Los autores agradecen el apoyo de la Universidad Nacional de Colombia, sede Bogotá y las excelentes condiciones de trabajo que se ofrecieron durante el ECOGyT, al comité organizador del evento y a las instituciones Yale University y Purdue University por el apoyo financiero. La segunda autora y el tercer autor agradecen el apoyo de NSF Grant DMS-2204365y  NSF Grant DMS-2405405, respectivamente.

\section{Equivalencias homotópicas, equivalencias débiles y grupos de homotopía}

En teoría de homotopía estudiamos \textit{espacios topológicos} salvo \textit{equivalencia homotópica} usando técnicas \textit{combinatorias} y \textit{algebraicas}. Para una discusión mas detallada de los resultados a continuación recomendamos los textos \cite{maunder, hatcher, gray, concise}.
\subsection{Homotopía}
\begin{defn} 
Sean $X$ y $Y$ dos espacios topológicos.
\begin{enumerate}
\item  Sean $f\colon X \to Y$ y $f'\colon X \to Y$ dos mapeos continuos. El mapeo $f$ es \textit{homotópico} a $f'$, denotado por $f \simeq f'$, si existe un mapeo continuo 
\[ h \colon X \times [0,1] \to Y\]
tal que $h(x,0)=f(x)$ y $h(x,1)=f'(x)$  para todo $x \in X$. Si esto se cumple, diremos que $h$ es una  \textit{homotopía de $f$ a $f'$}.
\item $X$ y $Y$ son \textit{homotópicamente equivalentes} si existen mapeos continuos $f \colon X \to Y$ y $g \colon Y \to X$ tales que $f \circ g \simeq \text{id}_Y$ y $g \circ f \simeq \text{id}_X$.
\end{enumerate}
\end{defn}

Intuitivamente, una homotopía $h$ puede ser vista como una familia 
\[\left\{h_t \colon X \to Y\right\}_{t \in [0,1]}\]
de mapeos continuos parametrizada por el intervalo $[0,1]$ que deforman continuamente $f=h_0$ en $f'=h_1$.
Notemos que la noción de homotopía define una relación de equivalencia en el conjunto de mapeos continuos entre los espacios $X$ y $Y$.
\begin{itemize}
    \item Reflexividad: Si $f:X\to Y$ es un mapeo continuo, entonces 
  $h(x,t)=f(x)$ es una homotopía de $f$ a $f$. 
    \item Simetría: Sean $f:X\to Y$ y $g:X\to Y$ dos funciones continuas. Suponga que $f\simeq g$, esto es, existe $h\colon X\times [0,1]\to Y$ continua, tal que $h(x,0)=f(x)$ y $h(x,1)=g(x)$. El mapeo $h'\colon X\times [0,1]\to Y$ definido por
    \[h'(x,t)=h(x,1-t)\]
     es una homotopía entre $g$ y $f$.
    \item Transitividad: Ejercicio.
    
\end{itemize}
Denotamos por $[X,Y]$ las clases de equivalencia del conjunto de mapeos continuos entre $X$ y $Y$ bajo la relación de homotopía. Es decir
\[ 
[X,Y]=\{f \colon X \to Y | f \text{ es un mapeo continuo  }\}/ \simeq 
\]

\begin{ex} \label{ejemplos1}
    \begin{enumerate} 
    \item Los espacios homotópicamente equivalentes a un punto se dicen \textit{contráctiles}. El disco 
    \[D^n=\{ x \in \mathbb{R}^n : ||x|| \leq 1\}\] 
    es contráctil, pero la esfera 
    \[S^{n-1}= \{ x \in \mathbb{R}^n : ||x|| =1\} \subset D^n\]
     no.
    Una de las motivaciones para la construcción de invariantes algebraicos es demostrar rigurosamente enunciados como este último. 
    
    \item Cualquier mapeo \[f \colon S^n \to S^n\] que no sea sobreyectivo es homotópico al mapeo constante \[c \colon S^n \to S^n\] que colapsa todo $S^n$ a un punto $* \in S^n$.
    
    \item Existe una biyección $[S^1, S^1]\cong \mathbb{Z}$. En general, para una espacio $X$ conexo por caminos existe una biyección entre $[S^1,X]$ y las clases de conjugación del grupo fundamental de $X$ (ver \cref{grphom}). 
    

    \item \textbf{Propiedad de levantamiento de homotopía para la inclusión $S^{n-1} \hookrightarrow D^n$.}\label{help} Dados un espacio $X$, una homotopía \[\{h_{t} \colon S^{n-1} \to X\}_{t\in [0,1]}\] y un mapeo $f \colon D^n \to X$, tal que $f|_{S^{n-1}}=h_0$, existe una homotopía $\{\widetilde{h}_t \colon D^n \to X\}$ extendiendo $\{h_t\}_{t\in [0,1]},$ i.e. satisfaciendo \[\widetilde{h}_t |_{S^{n-1}} = h_t\text{ para todo }t \in [0,1]\] con $\widetilde{h}_0=f$. Este es un enunciado fundamental en la formalización de la teoría de homotopía de espacios. Invitamos al lector a demostrarlo siguiendo la siguiente imagen.
    \[\includegraphics[scale=0.18]{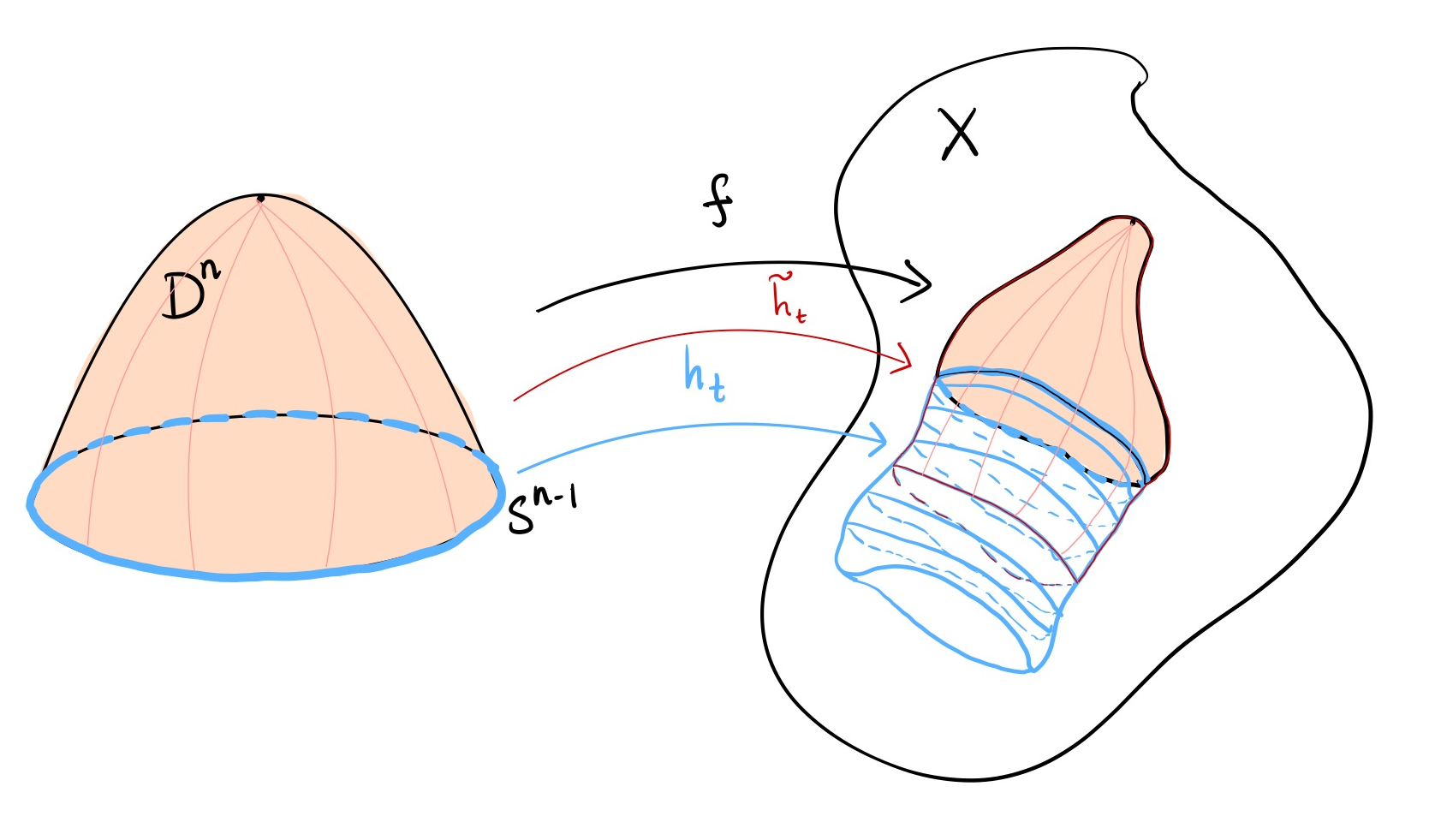}\]
    \item Claramente existen espacios que son homotópicamente equivalentes pero no homeomorfos (todo disco es homotópicamente equivalente, pero no homeomorfo, a un punto). Una tarea más interesante es analizar geométricamente la noción de equivalencia homotópica entre variedades de la misma dimensión. Por ejemplo, los \textit{espacios de Lens} $L(7,1)$ y $L(7,2)$ son dos variedades $3$-dimensionales que son homotópicamente equivalentes pero no homeomorfas, vea el ejercicio 2 de la sección 3.E y el ejercicio 29 de la sección 4.2 en \cite{hatcher}.  
      \end{enumerate}
\end{ex}

\subsection{Grupos de homotopía y CW-complejos}
Las equivalencias homotópicas entre espacios construidos pegando discos de diferentes dimensiones son detectadas por los \textit{grupos de homotopía}. 

\subsubsection{Grupos de homotopía}\label{grphom} Fijemos el punto base $*=(1,0, \ldots, 0) \in S^n$. Sea $X$ un espacio topológico y $b \in X$. Consideremos el conjunto de \text{mapeos puntuados}
\[ \text{Maps}((S^n,*), (X,b))=\{ \alpha \colon S^n \to X | \alpha \text{ es continuo y } \alpha(*)=b\}.\]
Para $\alpha,\alpha' \in \text{Maps}((S^n,*), (X,b))$, una homotopía $h \colon S^n \times [0,1] \to X$ entre $\alpha$ y $\alpha'$ es \textit{puntuada} si $h(*,t)=b$ para todo $t \in [0,1]$. Denotemos por $[(S^n,*), (X,b)] $ el conjunto de clases equivalencia de mapeos puntuados bajo la relación de homotopía puntuada. Cuando $n \geq 1$, podemos definir una estructura de grupo en $[(S^n,*), (X,b)]$ de la siguiente manera:

Sean $\alpha \colon S^n \to X$ y $\alpha' \colon S^n \to X$ representantes de dos clases en 
\[ [(S^n,*), (X,b)].\] Defina
\[\alpha \star \alpha' \colon (S^n,*) \to (X,b)\] 
como la composición
\[ S^n \to S^n \vee S^n \xrightarrow{\alpha \vee \alpha'}X, \]
donde el primer mapeo colapsa la copia de $S^{n-1}$ en $S^n$ dada por los puntos cuya última coordenada es $0$ y el segundo aplica $\alpha$ a la ``primera" copia de $S^n$ en $S^n \vee S^n$ y $\alpha'$ a la ``segunda". Este segundo mapeo lo aplicamos considerando el punto de identificación en $S^n \vee S^n$ como el punto base de ambas copias de $S^n$. 

\[\includegraphics[scale=0.15]{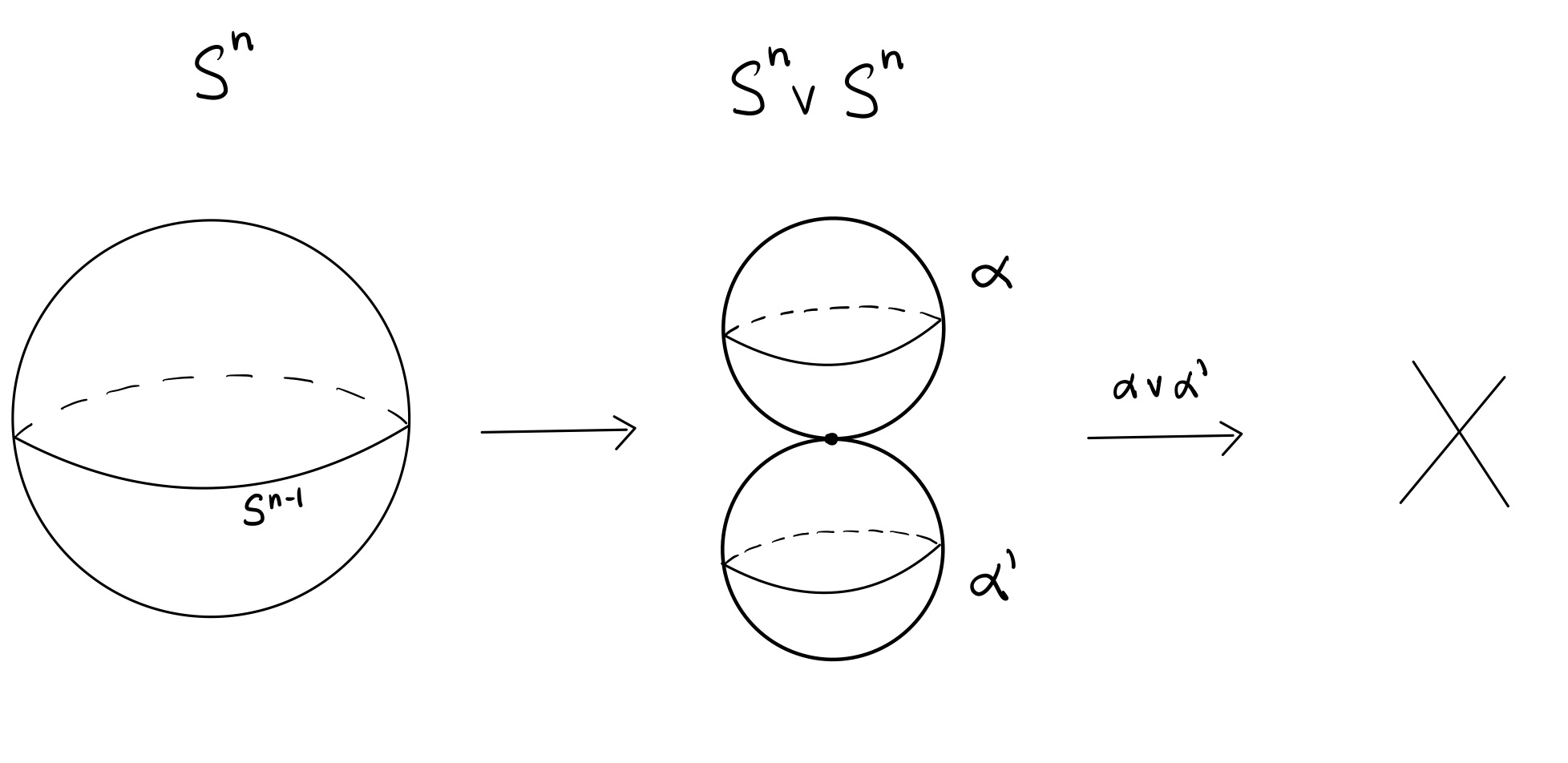}\]

\begin{propn}
Para todo entero $n \geq 1$, la operación binaria $\star$ está bien definida en el conjunto $[(S^n,*), (X,b)]$ y define una estructura de grupo.
\end{propn}

Para $n \geq 1$, el grupo $([(S^n,*), (X,b)], \star)$ es llamado $n$-\textit{grupo de homotopía} de $(X,b)$ y se denota por $\pi_n(X,b)$. El grupo $\pi_1(X,b)$ (que consiste en clases de homotopía puntuada de caminos que comienzan y terminan en $b$, con la operación binaria inducida por concatenación de caminos) se conoce como el \textit{grupo fundamental de $X$ en $b$} y fue introducido por Poincaré en el 1895 en el artículo \textit{Analysis Situs}. Este puede ser un grupo arbitrario, es decir, para cualquier grupo $G$ existe un CW-complejo puntuado $(X,b)$ tal que $\pi_1(X,b)\cong G$. Para aprender más sobre el grupo fundamental (y su íntima relación con espacios recubridores) sugerimos ver \cite{angie, massey}. 

Cuando $n\geq 2$ todos los grupos de homotopía $\pi_n(X,b)$ son grupos abelianos. Cuando $n=0$, $\pi_0(X)=[(S^0,*), (X,b)]$ es el conjunto de componentes conexas por caminos de $X$, en general no tiene estructura de grupo. 

Un mapeo puntuado $f \colon (X,b) \to (Y,c)$ induce un homomorfismo de grupos \[\pi_n(f) \colon \pi_n(X,b) \to \pi_n(Y,c)\] definiendo \[ \pi_n(f)[\alpha]= [ f \circ \alpha].\] Además, se tiene que $\pi_n(g \circ f) = \pi_n(g) \circ \pi_n(f)$ y $\pi_n(\text{id}_X)=\text{id}_{\pi_n(X,b)}$. 
\begin{propn} Si $f \colon (X,b) \to (Y,c)$ y $f' \colon (X,b) \to (Y,c)$ son mapeos homotópicos por medio de una homotopía puntuada entonces $\pi_n(f)=\pi_n(f')$ para todo $n \geq 0$.
\end{propn}
\begin{defn}
    Un mapeo continuo $f \colon X \to Y$ es una \textit{equivalencia débil de espacios topológicos} si \[\pi_n(f) \colon \pi_n(X,b) \to \pi_n(Y,f(b))\] es un isomorfismo para todo para todo $n \geq 0$ y $ b \in X$. 
\end{defn}

En general, es muy difícil calcular grupos de homotopía explícitamente. Todavía no sabemos todos los grupos de homotopía de las esferas y calcularlos sigue siendo uno de los problemas principales de la topología algebraica en el cual se han obtenido nuevos resultados recientemente. Algunos resultados básicos sobre grupos de homotopía de esferas son los siguientes. 

\begin{thm} \label{pi_n(S^n)}
Para todo $n \geq 1$ existen isomorfismos de grupos $\pi_n(S^n,*)\cong \mathbb{Z}$ y $\pi_k(S^n,*) \cong 0$ cuando $0< k <n$. 
\end{thm}

El siguiente resultado muestra que los grupos de homotopía de esferas exhiben un fenónmeno de \textit{estabilidad}.

\begin{thm} Para todo $n\geq 0$ y $k \geq 1$ existe un homomorfismo natural de grupos
\[\pi_{k}(S^n,*) \to  \pi_{k+1}(S^{n+1},*) \]
que es un isomorfismo para $k<2n-1$.
\end{thm}

\subsubsection{CW-complejos}
Un \textit{CW-complejo} es un espacio topológico $X$ construido de la siguiente manera:
\begin{enumerate}
    \item Comenzamos con un conjunto discreto $X^0$. Este conjunto lo pensamos como los ``vértices" del espacio que vamos a construir.
    \item Construimos inductivamente un espacio $X^n$ pegando un conjunto de $n$-discos $\{ D^n_{\alpha} \}$ a $X^{n-1}$ a través de mapeos continuos \[\phi_{\alpha} \colon S^{n-1} \cong \partial D^n_{\alpha} \to X^{n-1}\] usando la topología de cociente. Esto es, \[X^n = \big(X^{n-1} \bigsqcup_{\alpha} D^n_{\alpha}\big)/\sim \]
donde $x \sim \phi_{\alpha}(x)$ para todo $x \in \partial D^n_{\alpha}$. 
\item $X= \bigcup_{n=0}^{\infty} X^n$ con la siguiente topología: $U \subseteq X$ es un conjunto abierto si y solo si $U \cap X^n$ es un conjunto abierto para todo $n$.
\end{enumerate}

Un espacio topológico $X$ puede ser un CW-complejo de muchas maneras distintas. Cuando especificamos una descomposición particular dada por la información de una sucesión anidada de espacios \[X^0 \subseteq X^1 \subseteq X^2 \subseteq\cdots \subseteq X^n \subseteq\cdots\]
y una colección de mapeos continuos \[\{\varphi_{\alpha} \colon S^{n-1} \to X^{n-1}\}_{\alpha \in \mathcal{I}_n}\text{ para todo }n \geq 1\] que satisfacen los axiomas de la definición, la llamamos una \textit{estructura de CW-complejo para $X$}.

\begin{ex} Algunos ejemplos de espacios que admiten una estructura de CW-complejo son: grafos, la $n$-esfera $S^n$, el $n$-espacio proyectivo real $\mathbb{R}P^n$, el $n$-espacio proyectivo complejo $\mathbb{C}P^n$, el $n$-Grassmaniano (real o complejo) $\text{Gr}_n$, la realización geométrica de un complejo simplicial, variedades suaves compactas y  la realización geométrica de un conjunto simplicial (como veremos más adelante). 

Los siguientes son resultados más profundos. Toda variedad topológica cerrada de dimensión $n$, $n \neq 4$, admite una estructura de CW-complejo. La pregunta sobre si una variedad topológica cerrada de dimensión $4$ admite una estructura de CW-complejo sigue abierta. 
Sin embargo, no es tan difícil demostrar que toda variedad topológica compacta es \textit{homotópicamente equivalente} a un CW-complejo \cite[Corollary A.12]{hatcher}.

Si $X$ y $Y$ son CW-complejos tales que $X$ se obtiene pegando un conjunto finito de discos, y $Y$ pegando un conjunto contable de discos, entonces el espacio de mapeos $\text{Maps}(X,Y)$ (con la topología compacto-abierta) es homotópicamente equivalente a un CW-complejo. En particular, el espacio $\text{Maps}(S^k,S^n)$ es \textit{homotópicamente equivalente} a un CW-complejo \cite{milnorCW}.

\end{ex}

\begin{ex} Existen espacios topológicos que no admiten una estructura de CW-complejo. Un ejemplo de esto está dado por el \textit{círculo de Warsaw}: 
\[ W=\left\{\left(x,\text{sin}\left(\frac{1}{x}\right)\right): 0 < x \leq 1\right\} \cup \{(0,y) : -1\leq y \leq 1 \} \cup C, \] 
donde $C$ es un camino de $(0,-1)$ a $(1,\text{sin}(1))$. Note que aunque $W$ es conexo por caminos, no es localmente conexo por caminos y todos los CW-complejos lo son. 
\end{ex}

Las equivalencias homotópicas débiles son detectadas por conjuntos de clases de homotopía de mapeos desde todos los CW-complejos en el siguiente sentido.

\begin{thm}[Whitehead] \label{Whitehead}
Un mapeo continuo $f \colon X \to Y$ entre espacios topológicos es una equivalencia débil si y solo si para todo CW-complejo $Z$ el mapeo
\[f_*\colon [Z,X]\to [Z,Y] \] definido por
\[ f_* \colon (\phi\colon Z \to X) \mapsto (f \circ \phi \colon Z \to Y) \]
es biyectivo.
\end{thm}
Idea de la demostración: aplicar varias veces que la inclusión $S^n \hookrightarrow D^n$ satisface propiedades de \textit{extensión y levantamiento} de homotopías; vea el ejemplo \ref{help} y el Teorema ``HELP" en \cite{concise}.

\begin{cor}
Sean $X$ y $Y$ CW-complejos. Un mapeo continuo $f \colon X \to Y$ es una equivalencia débil si y solo si es una equivalencia homotópica. 
\end{cor}
\noindent \textit{Demostración.} Sea $f$ una equivalencia débil. Entonces, por el teorema anterior, el mapeo $[Y,X] \to [Y,Y]$ es sobreyectivo, por lo que existe algún $g \colon Y \to X$ tal que $f \circ g \simeq \text{id}_Y$. Esto  implica que $f \circ g \circ f\simeq f$. Pero también sabemos por el teorema anterior que el mapeo $[X,X] \to [X,Y]$ dado por $\phi \mapsto f \circ \phi $ es inyectivo, así que $g \circ f \simeq \text{id}_X$. La dirección conversa sigue del hecho de que dos mapeos homotópicos inducen el mismo mapeo en grupos de homotopía (teniendo cuidado con los puntos bases). \qed

\subsection{Motivación para estudiar teoría de homotopía}
 
La teoría de homotopía ha sido útil en otras áreas de la matemática. Algunos ejemplos son los siguientes. 

\begin{enumerate}
   \item Muchos problemas geométricos se reducen a calcular $[X,Y]$ para dos espacios $X$ y $Y$ que codifican o parametrizan información del problema en consideración y estructuras algebraicas en este conjunto. Por ejemplo, consideremos un espacio $X$ que se comporte de buena manera (paracompacto). Existe una biyección natural entre el conjunto de clases de isomorfismo de fibrados vectoriales de rango $n$ sobre $X$ y el conjunto $[X, \text{Gr}_n]$, donde el espacio $\text{Gr}_n$, llamado el $n$-\textit{Grasmanniano} y está definido como el conjunto
   \[ \text{Gr}_n=\{\text{subespacios $n$-dimensionles de } \mathbb{R}^{\infty} \}\] 
   equipado con una topología apropiada \cite{mitchell}.

   \item Problemas de clasificación de estructuras geométricas en  variedades salvo alguna noción equivalencia a menudo se reducen a calcular grupos de homotopía. Por ejemplo, el \textit{isomorfismo de Thom y Pontrjagin} nos dice que variedades suaves con una ``estructura tangencial" consideradas salvo \textit{cobordismo} están completamente determinadas por los grupos de homotopía de cierto espacio universal asociado a la estructura tangencial. \cite{diffviewpoint}

   \item Usando \textit{K-teoría algebraica} se pueden reducir algunos problemas difíciles y profundos de álgebra y teoría de números a calcular grupos de homotopía (usualmente difíciles de calcular explícitamente). 
   \item Las construcciones básicas de la teoría de homotopía de espacios se pueden abstraer usando teoría de categorías y convertirse en un formalismo aplicable a otras áreas de las matemáticas en donde la noción de isomorfismo se reemplaza por una mas flexible. Algunas áreas donde el formalismo abstracto de la teoría de homotopía se ha utilizado efectivamente son: álgebra homológica, física matemática, geometría simpléctica, teoría de representaciones, geometría algebraica, hasta ciencia de computación.

\end{enumerate}

\section{Conjuntos simpliciales}
La teoría de conjuntos simpliciales provee un modelo combinatorio para la teoría de homotopía de espacios topológicos. El concepto de \textit{conjunto simplicial} es una abstracción del concepto de \textit{complejo simplicial.} La idea central es describir un espacio (salvo equivalencia débil) con una ``receta combinatoria" de como construirlo. Algunas diferencias entre los conceptos de complejo simplicial y conjunto simplicial son 1) los símplices de un conjunto simplicial están coherentemente orientados, 2) los símplices de un conjunto símplicial no tienen que estar determinados por sus vértices y 3) los conjuntos simpliciales incluyen la información de símplices ``degenerados". En este documento solo discutiremos conjuntos simpliciales. Para una introducción al formalismo de conjuntos simpliciales con muchas ilustraciones y más detalles recomendamos \cite{friedman}. 

Los componentes básicos para definir los conjuntos simpliciales son los ordinales finitos y funciones no decrecientes entre ellos. 

\begin{notation}
Dado un entero $n\geq 0$, $[n]$ denota el conjunto ordenado $\{0<1<\dots < n\}.$ Para $m,n\geq 0$, $\mapdel([m],[n])$ denota el conjunto de funciones no decrecientes $[m] \to [n]$, i.e.
    
\[ \mapdel([m],[n]) = \{ f \colon [n] \to [m]
 | \text{ si $i \leq j$ entonces $f(i) \leq f(j)$} \}.\] 
\end{notation}

\begin{defn}
Un \emph{conjunto simplicial} $X$ consiste de
\begin{itemize}
    \item para cada $n\geq 0$, un conjunto $X_n$;
    \item para cada función no decreciente $f\colon [m] \to [n]$, una función
    \[f^*\colon X_n \to X_m,\]
\end{itemize}
tales que
\[\mathrm{id}^*=\mathrm{id} \qquad \text{y} \qquad (g\circ f)^* = f^*\circ g^*.\]
Otra notación común para $f^*$ es $X(f)$.
\end{defn}

\begin{ex}
 Dado $p\geq 0$, tenemos el conjunto simplicial $\Delta^p$, llamado el \emph{símplice estándar}, definido como
 \[(\Delta^p)_n=\mapdel([n],[p]).\]
 Para $f\colon [m] \to [n]$, el mapeo $f^* \colon (\Delta^p)_n \to (\Delta^p)_m$ es dado por composición.
\end{ex}

Los conjuntos simpliciales se pueden definir en términos de una colección más pequeña de mapeos.

\begin{notation}
    Dado $n$ y $0 \leq i \leq n$, consideremos la función 
    \[d^i \colon [n-1] \to [n]\]
    definida como
    \[d^i(x)=\begin{cases}
        x &\text{si } x<i\\
        x+1 &\text{si } x \geq i,
    \end{cases}\]
 es decir la función que salta $i$. Estas funciones se llaman \emph{cocaras}. Similarmente, consideremos la función 
 \[s^i \colon [n+1] \to [n]\]
    definida como
    \[s^i(x)=\begin{cases}
        x &\text{si } x\leq i\\
        x-1 &\text{si } x > i,
    \end{cases}\]
 es decir la función que repite $i$. Estas funciones se llaman \emph{codegeneraciones}.
\end{notation}

Esta colección de funciones genera todas las funciones no decrecientes $[m] \to [n]$, en el sentido del siguiente resultado. Estas funciones satisfacen una familia de identidades, llamadas las identidades cosimpliciales.

\begin{lemma}\label{cocaras}
    Toda función no decreciente $f\colon [m] \to [n]$ se puede escribir como el compuesto de cocaras y codegeneraciones. Las cocaras y codegeneraciones satisfacen las siguientes identidades:
    \begin{align*}
        d^jd^i&=d^id^{j-1} \text{ si } i<j\\
        s^jd^i&=d^is^{j-1} \text{ si } i<j\\
        s^jd^j&=\mathrm{id}=s^jd^{j+1}\\
        s^jd^i&=d^{i-1}s^j \text{ si } i >j+1\\
        s^js^i&=s^is^{j+1} \text{ si } i\leq j.
    \end{align*}
    
\end{lemma}

Este resultado implica directamente el siguiente corolario.

\begin{cor}
    Un conjunto simplicial $X$ está completamente determinado por la colección de conjuntos $X_n$ para $n\geq 0$, y para todo $0\leq i \leq n$, funciones
    \begin{align*}
        d_i & \colon X_n \to X_{n-1},\\
        s_i & \colon X_n \to X_{n+1},
    \end{align*}
    las cuales deben satisfacer las identidades simpliciales:
    \begin{align*}
        d_id_j&=d_{j-1}d_i \text{ si } i<j\\
        d_is_j&=s_{j-1}d_i \text{ si } i<j\\
        d_js_j&=\mathrm{id}=d_{j+1}s_j\\
        d_is_j&=s_jd_{i-1} \text{ si } i >j+1\\
        s_is_j&=s_{j+1}s_i \text{ si } i\leq j.
    \end{align*}
\end{cor}
Si $X$ es un conjunto simplicial, un elemento $x \in X_n$ se dice ser \textit{degenerado} si $x=s_i(y)$ para algún $0 \leq i \leq n-1$ y $y \in X_{n-1}$. De lo contrario, $x$ se dice ser \textit{no degenerado}. Por ejemplo, una cuenta sencilla revela que el número de elementos no degenerados en $(\Delta^p)_n$ es $\binom{p+1}{n}$ cuando $0 \leq n \leq p$  y $0$ cuando $n>p$.

\begin{ex} \label{BG} Dado un grupo $G$ con identidad $e$, definimos el \textit{espacio clasificante} $BG$ como el conjunto simplicial con $(BG)_n=G^n$, y caras y degeneraciones:
    \[d_i(g_1,\dots,g_n)=\begin{cases}
        (g_2,\dots,g_n) & \text{si } i=0;\\
        (g_1,\dots, g_{i+1} \cdot g_{i},\dots,g_n) & \text{si } 0<i<n;\\
        (g_1,\dots,g_{n-1}) & \text{si } i=n,
    \end{cases}\]
    \[s_i(g_1,\dots,g_n)=\begin{cases}
        (e,g_1,\dots,g_n) & \text{si } i=0;\\
        (g_1,\dots, g_i,e, g_{i+1},\dots,g_n) & \text{si } 0<i<n;\\
        (g_1,\dots,g_n,e) & \text{si } i=n.
    \end{cases}\]
De hecho, esta construcción se puede generalizar para de la siguiente manera. Dada una categoría $\mathsf{C}$ (ver \cref{cats}), definimos el \textit{nervio de $\mathsf{C}$} como el siguiente conjunto simplicial. Los $n$-símplices están dados por 
\[ (N \mathsf{C})_n= \{ (f_1,\dots,f_n) | f_i \in \mathsf{C}(\mathsf{s}(f_i), \mathsf{t}(f_i)) \text{ y } \mathsf{t}(f_i)=\mathsf{s}(f_{i+1})\},\]
donde $\mathsf{s}(f_i)$ y $\mathsf{t}(f_i)$ denotan el dominio y el rango de un morfismo $f_i$ en $\mathsf{C}$, respectivamente. En otras palabras, el conjunto de $n$-símplices de $N \mathsf{C}$ está dado por todas las secuencias componibles 
\[ x_0 \xrightarrow{f_1} x_1 \xrightarrow{f_2}  \cdots \xrightarrow{f_{n-1}} x_{n-1} \xrightarrow{f_n} x_n \] de $n$ morfismos en $\mathsf{C}$. Las caras y degeneraciones están definidas por:
    \[d_i(f_1,\dots,f_n)=\begin{cases}
        (f_2,\dots,f_n) & \text{si } i=0;\\
        (f_1,\dots, f_{i+1} \circ f_{i},\dots,f_n) & \text{si } 0<i<n;\\
        (f_1,\dots,f_{n-1}) & \text{si } i=n,
    \end{cases}\]
    \[s_i(f_1,\dots,f_n)=\begin{cases}
        (\text{id}_{\mathsf{s}(f_1)},f_1,\dots,f_n) & \text{si } i=0;\\
        (f_1,\dots, f_i,\text{id}_{\mathsf{t}(f_i)}, f_{i+1},\dots,f_n) & \text{si } 0<i<n;\\
        (f_1,\dots,f_n,\text{id}_{\mathsf{t}(f_n)}) & \text{si } i=n.
    \end{cases}\]
    Note que si consideramos un grupo $G$ como una categoría con un objeto entonces $BG= NG$. 
\end{ex}

\begin{ex}\label{sing}
Para todo espacio topológico $X$ podemos construir un conjunto simplicial llamado el \textit{complejo singular de $X$}, denotado por $\text{Sing}(X)$, definido de la siguiente manera. Primero consideremos el espacio topológico
\[ \mathbf{\Delta}^n = \{ (x_0,\ldots,x_n) \in \mathbb{R}^{n+1} | x_0 + \cdots +x_n=1, x_i \geq 0\} \]
llamado el $n$-\textit{símplice topológico}. Para todo $i=0,\ldots,n$ tenemos mapeos continuos (de hecho, lineales)
\[ d^i \colon \mathbf{\Delta}^{n-1} \to \mathbf{\Delta}^n\]
y 
\[s^i \colon \mathbf{\Delta}^{n+1} \to \mathbf{\Delta}^n\]
definidos por
\[ d^i(x_0,\ldots,x_n)= (x_0,\ldots,x_{i-1},0,x_{i+1}, \ldots, x_n)\]
y
\[ s^i(x_0, \dots, x_{n+1})= (x_0, \ldots, x_{i-1}, x_i + x_{i+1}, x_{i+2}, \ldots, x_n),\]
respectivamente. Los mapeos $d^i$ y $s^i$ satisfacen las identidades del Lema \ref{cocaras}. 
El conjunto simplicial $\text{Sing}(X)$ esta dado por
\[\text{Sing}(X)_n = \{ \sigma \colon \mathbf{\Delta}^n \to X | \sigma \text{ es un mapeo continuo } \} \]
con caras y degeneraciones
\[ d_i \colon \text{Sing}(X)_n\to \text{Sing}(X)_{n-1}\]
\[ s_i \colon \text{Sing}(X)_n \to \text{Sing}(X)_{n+1} \]
dados por
\[ d_i(\sigma)= \sigma \circ d^{i} \]
y
\[ s_i (\sigma) = \sigma \circ s^i,\]
respectivamente.

\end{ex}
\begin{ex}\label{realizacion}
Para todo conjunto simplicial $S$ podemos construir un espacio topológico llamado la \textit{realización geométrica de $S$}, denotado por $|S|$, definido de la siguiente manera:

\[ |S|= \coprod \big(S_n \times \mathbf{\Delta}^n \big) / \sim \]
donde $\sim$ es la relación de equivalencia generada por $(x,d^i(t)) \sim (d_i(x),t)$ para todo $(x,t) \in S_n \times \mathbf{\Delta}^{n-1}$ y $(y,s^i(u)) \sim  (s_i(y), u)$ para todo $(y,u) \in S_n \times \mathbf{\Delta}^{n+1}$. En otras palabras $|S|$ es el espacio topológico obtenido pegando una copia de $\mathbf{\Delta}^n$ para cada elemento de $S_n$ de acuerdo a las caras y degeneraciones de $S$.  
\end{ex}

\begin{ex} Sean $S$ y $S'$ conjuntos simpliciales con caras y degeneraciones denotadas por $d^S_i$, $s^S_i$ y $d^{S'}_i$, $s^{S'}_i$, respectivamente. El \textit{producto cartesiano de $S$ y $S'$}, denotado por $S \times S'$, es el conjunto simplicial dado por
\[(S \times S')_n = S_n \times S'_n\] con caras y degeneraciones dadas por
\[d^{S \times S'}_i = d^S_i \times d^{S'}_i\]
y
\[s^{S \times S'}_i= s^S_i \times s^{S'}_i.\]
\end{ex}

\begin{defn}
Si $X$ y $Y$ son conjuntos simpliciales, un \textit{mapeo simplicial} (o \textit{morfismo de conjuntos simpliciales}) $f \colon X \to Y$ consiste en un conjunto de funciones $\{f_n: X_n\to Y_n \}_{n \in \mathbb{Z}_{\geq 0}}$  tal que para todo morfismo $ \Theta:[m]\to[n]$, el siguiente diagrama 
    \[\begin{tikzcd}
X_n \arrow[rr, "f_n"] \arrow[dd, "\Theta^*=X(\Theta)"'] &  & Y_n \arrow[dd, "\Theta^*=Y(\Theta)"] \\
 &  &  \\
X_m \arrow[rr, "f_m"']                                  &  & Y_m                                 
\end{tikzcd}\]
conmuta; es decir $Y(\Theta) \circ f_n = f_m \circ X(\Theta)$.
\end{defn}

La teoría de conjuntos simpliciales se podría desarrollar sin la información de degeneraciones. Sin embargo, incluir esta información garantiza la existencia de un homeomorfismo natural de espacios topológicos
\[ |S \times S'| \cong |S| \times |S'|.\]

\section{Teoría de categorías}
La teoría de categorías es un lenguaje o marco teórico que nació en el contexto de topología algebraica con la motivación de establecer relaciones, analogías y equivalencias entre estructuras y construcciones provenientes de distintas áreas de las matemáticas. Para una introducción a la teoría de categoría recomendamos los textos \cite{maclane, emily}.
\begin{defn}\label{cats}
    Una \textit{categoría} $\mathsf{C}$ consiste de la siguiente información:
\begin{itemize}
    \item una colección $Ob(\mathsf{C})$
    cuyos elementos se llaman \textit{objetos},
    \item para todo par de objetos $x$ y $y$ en $Ob(\mathsf{C})$, una colección $\mathsf{C}(x,y)$ cuyos elementos se llaman \textit{morfismos},
    
    \item una \textit{regla de composición} que para todo $x,y,z \in Ob(\mathsf{C})$, $f \in \mathsf{C}(x,y)$ y $g \in \mathsf{C}(y,z)$, asigna un morfismo $g \circ f \in \mathsf{C}(x,z)$ y
    \item para todo objeto $x \in Ob(\mathsf{C})$ un morfismo $\text{id}_x \in \mathsf{C}(x,x)$ llamado la \textit{identidad en} $x$. 
\end{itemize}
Esta información se requiere estar sujeta a las siguientes propiedades:

\begin{itemize}
\item \textit{la regla de composición es asociativa}: para todo $f\in \mathsf{C}(x,y), g \in \mathsf{C}(y,z)$ y $h \in \mathsf{C}(z,w)$ se tiene que \[(h \circ g) \circ f= h \circ (g \circ f).\]

\item \textit{la regla de composición es unitaria}: para todo $f \in \mathsf{C}(x,y)$, se tiene que \[\text{id}_y \circ f=f=f \circ \text{id}_x.\]
\end{itemize}
En general, la colección de objetos una categoría no forman un conjunto, más bien, una clase propia en el sentido de teoría de conjuntos. Este detalle técnico no debería preocupar al lector por ahora, ya que nunca consideraremos todos los objetos de una categoría cuyos objetos no forman un conjunto como una totalidad para hacer alguna construcción. En la mayoría de nuestros ejemplos, para todo par de objetos $x$ y $y$ en una categoría $\mathsf{C}$, la colección de morfismos $\mathsf{C}(x,y)$ será un conjunto. Esto es un principio importante de la teoría de categorías al nivel que la estaremos discutiendo: en muchas ocasiones no nos interesa trabajar con todos los objetos a la vez pero sí con los morfismos entre cada pareja de objetos.

Un morfismo $f \in \mathsf{C}(x,y)$ es un \textit{isomorfismo} si existe $f^{-1} \in \mathsf{C}(y,x)$ tal que $f \circ f^{-1} = \text{id}_y$ y $f^{-1} \circ f= \text{id}_x$.

Aveces abusaremos un poco de la notación y escribiremos $x \in \mathsf{C}$ cuando $x$ es un objeto de $\mathsf{C}$.  

\end{defn}
Muchas clases de objetos matemáticos familiares forman una categoría como muestran los siguientes ejemplos.
\begin{enumerate}
    \item La categoría $\mathsf{Set} $ tiene conjuntos como objetos y funciones como morfismos.
   
    \item La categoría $\mathsf{Top}$ tiene espacios topológicos como objetos y mapeos continuos como morfismos.
    \item La categoría $\mathsf{hoTop}$ tiene espacios topológicos como objetos y clases de homotopía de mapeos continuos como morfismos, i.e. $\mathsf{hoTop}(X,Y)=[X,Y]$. 
    \item La categoría $\mathsf{Man}$ tiene variedades suaves como objetos y mapeos suaves como morfismos. 
    \item La categoría $\mathsf{Gp}$ tiene grupos como objetos y homomorfismos de grupo como morfismos. Las categorías $\mathsf{Mon}$, $R\text{-}\mathsf{Mod}$, y $F\text{-}\mathsf{Vect}$, $\mathsf{Rings}$ de monoides, módulos (derechos o izquierdos) sobre un anillo $R$, espacios vectoriales sobre un campo $F$ y anillos, respectivamente, se definen de manera similar; en todos los casos los morfismos son mapeos que preservan la estructura algebraica.
    
    \item La categoría $\Delta$ tiene como objetos 
    \[Ob(\Delta)= \{[n] | n\in \mathbb{Z}_{\geq 0} \}\] y funciones $f \colon [n] \to [m]$ no-decrecientes como morfismos.
     \item La categoría $\mathsf{sSet}$ tiene conjuntos simpliciales como objetos y mapeos simpliciales como morfismos. 
    \item Para todo monoide $M$ podemos definir una categoría $\mathsf{B}M$ con un solo objeto $*$ y morfismos $\mathsf{B}M(*,*)=M$. La regla de composición está determinada por la multiplicación de $M$ y la identidad en $*$ por la identidad de $M$. 

    \item \label{caminos} Para todo espacio topológico $X$ podemos definir una categoría $\mathsf{P}X$ de la siguiente manera. Los objetos de $\mathsf{P}X$ son todos los puntos de $X$, es decir, $Ob(\mathsf{P}X)$ es el conjunto subyacente de $X$. Para todo $a,b \in Ob(\mathsf{P}X)$ definimos un conjunto de morfismos como
    \begin{eqnarray*}
    \mathsf{P}X(a,b) = \{ (\gamma, r) | r \in [0, \infty), \gamma \colon [0,r] \to X \\ \text{ es un mapeo continuo y } \gamma(0)=a, \gamma(r)=b \}.
       \end{eqnarray*}
Es decir, $\mathsf{P}X(a,b)$ es el conjunto de caminos continuos de $a$ a $b$ parametrizados por un intervalo $[0,r]$ para algún parámetro $r$ (no fijo). La regla de composición está definida por concatenación de caminos sumando los parámetros correspondientes. La identidad en $a \in X$ es el camino constante $id_a \colon [0,0]=\{0\} \to X$, $id_a(0)=a$. Los conjuntos de morfismos $\mathsf{P}X(a,b)$ pueden ser considerado como espacios topológicos con la topología inducida por la topología compacto-abierta. Bajo esta topología, la composición de morfismos $\mathsf{P}X$ es un mapeo continuo. 
\item \label{grupoidefundamental} Sea $\pi_0(\mathsf{P}X(a,b))$ el conjunto de componentes conexas por caminos del espacio topológico $\mathsf{P}X(a,b)$ del ejemplo anterior. Para todo espacio topológico $X$ defina una categoría $\Pi(X)$ cuyos objetos son todos los puntos de $X$ y para todo $a,b \in X$, \[\Pi(X)(a,b)= \pi_0(\mathsf{P}X(a,b)).\] La composición está inducida por concatenación de caminos y la identidad en $a \in X$ es la clase de equivalencia de $id_a$ en  $\pi_0(\mathsf{P}X(a,a))$. Todo morfismo en $\Pi(X)$ es un isomorfismo. El inverso de un morfismo es representado por recorrer el camino correspondiente en la dirección opuesta. La categoría $\Pi(X)$ se llama el \textit{grupoide fundamental de $X$}. Note que podemos identificar $\Pi(X)(b,b)$ con el grupo fundamental $\pi_1(X,b)$.

\item Un conjunto parcialmente ordenado $(P, \leq)$ define una categoría con objetos los elementos de $P$ y para todo par $i,j \in P$ un morfismo único de $i$ a $j$ si y solo si $i\leq j$.

\item Un grafo dirigido $G$ define una categoría $PG$ con objetos los vértices de $G$ y, para todo par de vértices $v$ y $w$, $PG(v,w)$ es el conjunto de caminos dirigidos de $v$ a $w$. Los morfismos de identidad se añaden formalmente para cada vértice $v \in G$.  La regla de composición es dada por concatenación de caminos. 

    \item Para toda categoría $\mathsf{C}$ podemos definir una nueva categoría $\mathsf{C}^{op}$ con morfismos $\mathsf{C}^{op}(x,y)= \mathsf{C}(y,x)$ y regla de composición $f \circ^{op} g = g \circ f.$ En otras palabras, $\mathsf{C}^{op}$ es la categoría que obtenemos virando la dirección de todos los morfismos de $\mathsf{C}$. 
\end{enumerate}

\begin{defn}
    Sean $\mathsf{C}$ y $\mathsf{D}$  categorías. Un \textit{funtor} $F:\mathsf{C}\to \mathsf{D}$ consiste en la información de
    \begin{itemize}
        \item un objeto $F(x) \in Ob(\mathsf{D})$ para todo  objeto $x \in Ob(\mathsf{C})$ y
        \item un morfismo $F(f) \colon F(x) \to F(y)$ en $\mathsf{D}$ para todo morfismo $f \colon x \to y$ en $\mathsf{C}$
        \end{itemize}
tal que para todo par de morfismos $f \colon x \to y$ y $g \colon y \to z$ en $\mathsf{C}$ tenemos que
 \[F(g\circ f)= F(g) \circ F(f)\]
 y para todo objeto $x \in \mathsf{C}$ tenemos que
   \[ F(\text{id}_x)=id_{F(x)}. \]
\end{defn}
La colección de categorías se puede organizar como una categoría. Denotemos por $\mathsf{Cat}$ la categoría cuyos objetos son categorías y cuyos morfismos son funtores.

 Una \textit{categoría enriquecida en espacios topológicos} consiste en un categoría $\mathsf{C}$ con la estructura adicional de una topología en cada conjunto de morfismos $\mathsf{C}(x,y)$ tal que la regla de composición
    \[ \mathsf{C}(y,z) \times \mathsf{C}(x,y) \to \mathsf{C}(x,z) \]
    \[ (g,f) \mapsto g \circ f \]
    es un mapeo continuo. Las categorías enriquecidas forman una categoría $\mathsf{Cat}_{\mathsf{Top}}$ con morfismos los funtores $F \colon \mathsf{C} \to \mathsf{D}$ tal que el mapeo inducido en morfismos \[F \colon \mathsf{C}(x,y) \to \mathsf{D}(F(x), F(y))\] es continuo. Por ejemplo, para todo espacio topológico $X$, la construcción $\mathsf{P}(X)$ del ejemplo \ref{caminos} define una categoría enriquecida en espacios topológicos.

\begin{ex} \label{funtores}
Continuamos con una lista de ejemplos de funtores. Note que en todos estos ejemplos hay una manera ``natural" de definir cada funtor al nivel de morfismos. 
\begin{enumerate}
    \item Funtor de componentes conexos
    \[ \pi_0 \colon  \mathsf{Top} \to \mathsf{Set} \]
    y funtor de grupo de homotopía \[\pi_n:\mathsf{Top}_*\to  \mathsf{Gp} \] para cada $n\geq 1$
    \item Funtor de grupo de homología singular
    \[H_n \colon \mathsf{Top} \to \mathsf{\mathbb{Z}\text{-}\mathsf{Mod}}\] para cada $n \geq 0$ (el cual no hemos definido en este documento pero es un invariante importante de espacios topológicos; vea \cite{hatcher, concise, vick})
    \item La construcción $X \mapsto \mathsf{P}(X)$ del ejemplo \ref{caminos} define un funtor
   \[  \mathsf{P} \colon \mathsf{Top} \to \mathsf{Cat}_{\mathsf{Top}}.\]
   \item El grupoide fundamental $X \mapsto \Pi(X)$, construido en el ejemplo \ref{grupoidefundamental}, define un funtor
   \[ \Pi \colon \mathsf{Top} \to \mathsf{Cat}.\]

   \item \label{olvidogrupos} El funtor de olvido \[\mathcal{U}:  \mathsf{Gp} \to \mathsf{Set}\] ``olvida" las operación binaria de un grupo y solo ``recuerda" el conjunto subyacente.
   
   \item \label{olvidoespacios} El funtor de olvido \[\mathcal{U} \colon \mathsf{Top} \to \mathsf{Set}\] ``olvida" la topología de un espacio y solo ``recuerda" el conjunto subyacente. 
    
    \item \label{olvidocat} El funtor de olvido \[\mathcal{U} \colon \mathsf{Top} \to \mathsf{Set}\] induce otro funtor (abusando de la notación) \[\mathcal{U} \colon  \mathsf{Cat}_{\mathsf{Top}} \to \mathsf{Cat}\] que olvida la estructura de espacio topológico en los morfismos de una categoría enriquecida en espacios topológicos. 

        \item \label{libregrupos} El funtor grupo libre \[F:\mathsf{Set} \to  \mathsf{Gp}\] que envía un conjunto $S$ al grupo libre $F(S)$ generado por $S$. 
        
    \item \label{libreespacios} El funtor \[F \colon \mathsf{Set} \to \mathsf{Top}\] que envía un conjunto $S$ al espacio topológico dado por $S$ con la topología discreta (todos los subconjuntos de $S$ son abiertos). 

    \item \label{librecat} El funtor $F \colon \mathsf{Set} \to \mathsf{Top}$ induce un funtor \[ F \colon \mathsf{Cat} \to \mathsf{Cat}_{\mathsf{Top}}.\]

    \item Sea $X$ un espacio topológico y $\mathcal{O}_X$ la categoría cuyos objetos son conjuntos abiertos en $X$ y morfismos son inclusiones. Tenemos un funtor \[C \colon \mathcal{O}_X^{op} \to \mathsf{Rings}\] que asigna a cada abierto $U \subseteq X$ el anillo \[C(U,\mathbb{R})= \{ f \colon U \to \mathbb{R}| f \text{ es una función continua} \}\]
    y a cada inclusión $i \colon U \hookrightarrow V$ el mapeo de restricción $i^* \colon C(V,\mathbb{R}) \to C(U,\mathbb{R})$ correspondiente. 
\item El nervio de una categoría, discutido en el ejemplo \ref{BG}, define un funtor
\[N \colon \mathsf{Cat} \to \mathsf{sSet}.\]
Note que, si consideramos $[n]=\{0 \to 1 \to \cdots \to n\}$ como una categoría, tenemos que el conjunto de $n$-símplices del nervio está dado por
\[ N(\mathsf{C})_n= \mathsf{Cat}([n], \mathsf{C}).\]
    
    \item El complejo singular de un espacio, discutido en el ejemplo \ref{sing}, define un funtor
    \[ \text{Sing} \colon \mathsf{Top} \to \mathsf{sSet}.\]
    \item La realización geométrica de un conjunto simplicial, discutida en el ejemplo \ref{realizacion}, define un funtor
    \[ |-| \colon \mathsf{sSet} \to \mathsf{Top}.\]
       \end{enumerate}
\end{ex}

Notemos que un conjunto simplicial se puede interpretar como un funtor
\[X:\Delta^{op}\to \mathsf{Set}\]
con $X([n])=X_n$. Bajo esta perspectiva, los morfismos de conjuntos simpliciales son transformaciones naturales de funtores, una noción definida de la siguiente manera.

\begin{defn}
    Sean $\mathsf{C}$ y $\mathsf{D}$ categorías y 
    \[F,G: \mathsf{C}\to \mathsf{D}\]
   dos funtores. Una \textit{transformación natural} $\alpha$ de $F$ a $G$, denotada por
    \[\alpha: \mathsf{C} \Rightarrow \mathsf{D}, \]
    es una colección de morfismos 
        \[\alpha_x:F(x) \to G(x) \]
        para todo objeto $x$ en $\mathsf{C}$, tal que para todo morfismo $f:x\to y$ en $\mathsf{C}$, el siguiente diagrama
        \[\begin{tikzcd}
F(x) \arrow[rr, "\alpha_x"] \arrow[dd, "F(f)"'] &  & G(x) \arrow[dd, "G(f)"] \\
                                                &  &                         \\
F(y) \arrow[rr, "\alpha_y"']                     &  & G(y)                   
\end{tikzcd}\]
conmuta; es decir, $G(f) \circ \alpha_x= \alpha_y \circ F(f)$.
  
\end{defn}
Dadas dos categorías $\mathsf{C}$ y $\mathsf{D}$ podemos definir una nueva categoría $\mathsf{Fun}(\mathsf{C}, \mathsf{D})$ con funtores $\mathsf{C} \to \mathsf{D}$ como objetos y transformaciones naturales como morfismos. En particular, $\mathsf{sC}=\mathsf{Fun}(\Delta^{op}, \mathsf{C})$ se conoce como la categoría de \textit{objetos simpliciales en $\mathsf{C}$}.

 \begin{ex} Algunos ejemplos de transformaciones naturales son los siguientes:
 
 \begin{enumerate}
 \item Sean $F:\mathsf{Set} \to  \mathsf{Gp} $ y $\mathcal{U}:  \mathsf{Gp} \to \mathsf{Set}$ los funtores de los ejemplos \ref{libregrupos} y \ref{olvidogrupos}, respectivamente. Para todo conjunto $S$ definimos un morfismo $\eta_S:S\to \mathcal{U}F(S)$ por $\eta_S(s)=s$. Para todo grupo $G$, definimos un morfismo $\varepsilon_G \colon F \mathcal{U}(G) \to G$ por 
        $\varepsilon_G(g_1 \ldots g_n) = g_1 * \ldots * g_n$ donde $*$ denota la multiplicación de $G$. Las colecciones de morfismos $\eta_S$ y $\varepsilon_G$ definen transformaciones naturales
        \[\eta:\text{id}_{Set}\Rightarrow\mathcal{U}F\]
        y
        \[\varepsilon:F\mathcal{U}\Rightarrow \text{id}_{ \mathsf{Gp} },\]
        respectivamente. 
\item De manera similar al ejemplo anterior podemos definir transformaciones naturales $\text{id} \Rightarrow \mathcal{U} F$ y $\varepsilon \colon F\mathcal{U} \Rightarrow \text{id}$ para las parejas de funtores $(\mathcal{U}, F)$ de los ejemplos \ref{libreespacios} y \ref{olvidoespacios} y de los ejemplos \ref{librecat} y \ref{olvidocat}.
        \item Considere los funtores

\[  \mathsf{Top}_* \xrightarrow{\pi_n} \mathsf{Gp}\]
y
\[ \mathsf{Top}_* \xrightarrow{J} \mathsf{Top} \xrightarrow{H_n} \mathsf{Gp} \]
donde el funtor$J$ simplemente olvida el punto base de un espacio puntuado. Para todo espacio puntuado $X$, sea \[\mathsf{h}_X \colon \pi_n(X) \to H_n(JX)\] el morfismo de grupos determinado de la siguiente manera. Sea $\kappa_n \in H_n(S^n) \cong \mathbb{Z}$ un generador de la $n$-homología de la $n$-esfera. Para una clase $[\sigma \colon S^n \to X] \in \pi_n(X)$ defina $\mathsf{h}_X(\sigma) = H_n(\sigma)(\kappa_n) \in H_n(JX)$. La colección de morfismos $\mathsf{h}_X$ determina una transformación natural \[\mathsf{h} \colon \pi_n \Rightarrow H_nJ\] conocida como \textit{la transformación natural de Hurewicz}.

\item Siguiendo los ejemplos en \ref{funtores}, considere los funtores
\[ \mathsf{Top} \xrightarrow{\mathsf{P}} \mathsf{Cat}_{\mathsf{Top}}  \]
y
\[ \mathsf{Top}\xrightarrow{\Pi} \mathsf{Cat} \xrightarrow{F} \mathsf{Cat}_{\mathsf{Top}}.\]
Para todo espacio $X$ consideremos el funtor $\pi_X \colon \mathcal{U}(\mathsf{P}(X)) \to \Pi(X)$ que es identidad en objetos y a nivel de morfismos de una categoría enriquecida en espacios topológicos aplica el funtor de componentes conexos con la topología discreta. La colección de funtores $\pi_X$ determina una transformación natural \[ \pi \colon \mathsf{P} \Rightarrow F\Pi .\]

        \end{enumerate}
\end{ex}

\section{Equivalencia entre espacios topológicos y conjuntos simpliciales}

Para todo conjunto simplicial $S$ hemos definido un espacio topológico $|S|$ y para todo espacio topológico $X$ un conjunto simplicial $\Sing(X)$. Discutiremos las siguientes preguntas:
\begin{enumerate}
    \item ¿Qué tipo de espacio topológico es $|S|$ y qué tipo conjunto simplicial es $\Sing(X)$? ¿Cuáles son sus propiedades esenciales?
    \item ¿Cuánto``recuerda" $\Sing(X)$ de $X$ y $|S|$ de $S$? ¿En qué sentido $\mathsf{sSet}$ y $\mathsf{Top}$ son equivalentes?
    \item ¿En qué sentido un espacio $X$ o conjunto simplicial $S$ se parece a una categoría?
\end{enumerate}

Para una discusión más detallada sobre los resultados de esta sección vea \cite{friedman, goerss-jardine, simplicialobjects, curtis, kerodon}.
\subsection{Realización geométrica y complejo singular}

Recordemos que en los ejemplos \ref{realizacion} y \ref{sing} definimos los funtores de realización geométrica y complejo singular, respectivamente. Para todo conjunto simplicial $S$ el espacio topológico $|S|$ es construido pegando copias de símplices topológicos de manera similar a la que un CW-complejo se construye pegando copias de discos y por ende no es difícil concluir lo siguiente.

\begin{propn} Para todo conjunto simplicial $S$, $|S|$ es un CW-complejo.
\end{propn}

 Denotemos por $\iota_n$ el único símplice no-degenerado en $(\Delta^n)_n$. Este elemento coresponde a la identidad $\text{id}_{[n]} \colon [n] \to [n]$ bajo la definición \[(\Delta^n)_n =\Delta([n],[n]).\] 
De hecho, existe una biyección natural 
\begin{eqnarray}\label{yoneda}
S_n \cong \mathsf{sSet}(\Delta^n,S)
\end{eqnarray}
para todo $S \in \mathsf{sSet}$ que envía un elemento $\sigma \in S_n$ al morfismo de conjuntos simpliciales $\Delta^n \to S$ determinado por $\iota_n\mapsto \sigma$. Note que todos los símplices no-degenerados de $\Delta^n$ se pueden expresar como $d_{i_1} \ldots d_{i_k}\iota_n \in (\Delta^n)_{n-k}$ para algún $0 \leq k \leq n$. Estas expresiones no son únicas; las redundancias están dadas por las identidades simpliciales que satisfacen las composiciones de caras. Los símplices degenerados en $\Delta^n$ están determinados formalmente por los no-degenerados.

Para describir qué tipo de objeto es $\text{Sing}(X)$ es útil introducir una colección de conjuntos simpliciales llamados \textit{cuernos}. 
\begin{defn}
    Para todo $n \geq 0$, definimos el $(n,i)$-\textit{cuerno} para $0 \leq i \leq n$ como el subconjunto simplicial $\Lambda^n_i \subset \Delta^n$ determinado por todos los símplices no-degenerados de $(\Delta^n)_{n-1}$ excepto $d_i \iota_n$. 
\end{defn}

\begin{defn} \label{Kan}
Un \textit{complejo de Kan} $S$ es un conjunto simplicial con la siguiente propiedad:
para todo $n \geq 1$ y todo morfismo $f\colon \Lambda^n_i \to S$ para $0 \leq i \leq n$, existe una extensión $\tilde{f} \colon \Delta^n \to S$; es decir, existe un morfismo $\tilde{f} \colon \Delta^n \to S$ tal que $\tilde{f}|_{\Lambda^n_i}=f$. Llamamos a $\tilde{f}$ un \textit{relleno} para $f$. 
\end{defn}

Un complejo de Kan $K$ puede ser interpretado
como un tipo de categoría en donde la composición de morfismos está definida salvo ``homotopías coherentes" y en donde cada morfismo es invertible salvo homotopía. De hecho, los elementos de $K_0$ pueden ser pensados como objetos y un elemento $\alpha = (\alpha \colon \Delta^1 \to K) \in K_1$ (utilizando la identificación \ref{yoneda}) como un morfismo de $d_1\alpha \in K_0$ a $d_0\alpha \in K_0$. Dos morfismos componibles, es decir dos elementos $\alpha$ y  $\beta$ en $K_1$ tal que $d_0\alpha= d_1\beta$, determinan un mapeo simplicial \[\alpha \wedge \beta \colon \Lambda^2_1 \to K.\] Como $K$ es un complejo de Kan, existe un relleno 
\[f \colon \Delta^2 \to K\]
para $\alpha \wedge \beta$. Entonces \[ \Delta^1 \xrightarrow{d^1} \Delta^2 \xrightarrow{f} K \] determina un elemento de $K_1$ que puede ser considerado como la composición $\beta \circ \alpha$. Pero el relleno  $f$ no tiene por qué ser único, así que esta composición no está bien definida. Sin embargo, si $f' \colon \Delta^2 \to K$ otro relleno para $\alpha \wedge \beta$, entonces $f $ y $f'$ determinan un mapeo simplicial 
\[ h \colon \Lambda^3_1 \to K \]
tal que la restricción de $h$ a la cara $d_3\iota_3$ es $f$, la restricción de $h$ a la cara $d_2\iota_3$ es $f'$ y la restricción de $h$ a la cara $d_0\iota_3$ es $\Delta^2 \xrightarrow{s^1} \Delta^1 \xrightarrow{\beta} K$. Como $K$ es un complejo de Kan, existe un relleno \[\tilde{h} \colon \Delta^3 \to K\] para $h$. Entonces 
\[ \Delta^2 \xrightarrow{d^1} \Delta^3 \xrightarrow{\tilde{h}} K\]
determina un elemento de $K_2$ que puede ser considerado como una homotopía entre $f$ y $f'$. Pero esta homotopía no es única y dada otra homotopía, por la condición de Kan, existe una homotopía superior en $K_3$ relacionándolas y así sucesivamente. También note que cualquier elemento $\alpha \in K_1$ determina un mapeo simplicial
\[ i_{\alpha} \colon \Lambda^2_0 \to K\]
tal que la restricción de $i_{\alpha}$ a la cara $d_2\iota_2$ es $\alpha$ y la restricción de $i_{\alpha}$ a la cara $d_1\iota_2$ es
$s_0d_1\alpha \colon \Delta^1 \xrightarrow{s^0} \Delta^0 \xrightarrow{d^1} \Delta^1 \xrightarrow{\alpha} K$. Como $K$ es un complejo de Kan, existe un relleno 
$\tilde{i}_{\alpha} \colon \Delta^2 \to K$. Entonces
\[ \Delta^1 \xrightarrow{d^0} \Delta^2 \xrightarrow{\tilde{i}_{\alpha}} K\]
determina un elemento en $K_1$ que puede ser considerado como un ``inverso izquierdo" para $\alpha$ salvo la homotopía $\tilde{i}_{\alpha}$, si pensamos en $s_0d_1\alpha$ como el morfismo identidad en $d_1\alpha \in K_0$. Similarmente, podemos obtener un ``inverso derecho" para $\alpha$ salvo homotopía y ambos inversos están relacionados por una homotopía y estas homotopías están relacionadas por homotopías superiores y así sucesivamente. 

El ejemplo mas importante de un complejo de Kan es el complejo singular de un espacio topológico.

\begin{propn} \label{SingesKan} Para todo espacio topológico $X$, $\Sing(X)$ es un complejo de Kan.
\end{propn}

Otro ejemplo de un complejo de Kan es $BG$ para cualquier grupo $G$ y los rellenos en este caso son únicos. En general, $N\mathsf{C}$ no es un complejo de Kan para una categoría arbitraria $\mathsf{C}$. Sin embargo, $N\mathsf{C}$ satisface la condición de rellenos para todo $f \colon \Lambda^n_i \to N\mathsf{C}$ si $0 < i < n$ y los rellenos en este caso son únicos.

Los conjuntos simpliciales que satisfacen la condición de relleno de la definición \ref{Kan} pero solo para  $0 < i < n$ se llaman \textit{cuasi-categorías} \cite{lurie-what}. Una cuasi-categoría la podemos pensar como una categoría en la que la composición de morfismos está definida salvo homotopías coherentes. Note que sin poder rellenar los $(n,0)$- y $(n,n)$-cuernos no podemos concluir que los morfismos (o los $1$-símplices) son invertibles salvo homotopía. 

\subsection{Una equivalencia de \textit{teorías de homotopía}}
Las nociones de homotopía y grupos de homotopía también pueden ser definidas en el contexto simplicial.  
\begin{defn}
Sean $f_0,f_1 \colon S \to S'$ dos morfismos de conjuntos simpliciales. Una \textit{homotopía simplicial de $f_0$ a $f_1$} es un morfismo de conjnuntos simpliciales 
\[ h \colon S \times \Delta^1 \to S'\]
tal que la composición de morfismos
\[S \cong S \times \Delta^0 \xrightarrow{\text{id}_S \times d^i} S \times \Delta^1 \xrightarrow{h} S'\] es exactamente $f_{1-i} \colon S \to S'$, para $i=0,1$.

\end{defn}

\begin{propn}
Si $K$ es un complejo de Kan, la relación inducida en el conjunto $\mathsf{sSet}(S,K)$ por la noción de homotopía simplicial es una relación de equivalencia para todo conjunto simplicial $S$. 
\end{propn}

Al igual que en el caso de espacios, denotamos por $[S,K]$ el conjunto de clases de equivalencia de mapeos de conjuntos simpliciales salvo homotopía simplicial. 

También podemos definir combinatóricamente los grupos de homotopía para todo complejo de Kan. Sea $\partial \Delta^n$ el subconjunto simplicial de $\Delta^n$ que obtenemos removiendo $\iota_n \in (\Delta^n)_n$ (y sus degeneraciones correspondientes). Si $K$ es un complejo de Kan y $v \in K_0$, definimos $\pi_n(K,v)$ como el conjunto de clases de equivalencia $[\alpha]$ donde $\alpha \colon \Delta^n \to K$ es un mapeo simplicial que colapsa todos los $i$-símplces de $\partial \Delta^n \subset \Delta^n$ a $(s_0 \circ \ldots \circ s_0)(v) \in K_i$ (para $0 \leq i \leq n-1$) y la relación de equivalencia está dada por homotopíaas simpliciales $h \colon \Delta^n \times \Delta^1 \to K$ que de la misma manera colapsan todos los símplices de $\partial \Delta^n \times \Delta^1$ a la degeneración correspondiente de $v \in K_0$. La estructura de grupo se define de tal manera que existe un isomorfismo natural de grupos $\pi_n(K,v) \cong \pi_n(|K|,v)$, para más detalles vea \cite[\href{https://kerodon.net/tag/00VJ}{Tag 00VJ}]{kerodon}.

\begin{defn} \label{equivdebilsset}
Sean $S$ y $S'$ conjuntos simpliciales. Un mapeo $f \colon S\to S'$ es una \textit{equivalencia débil de conjuntos simpliciales} si para todo complejo de Kan $K$, el mapeo de conjuntos
\[ f^* \colon [S',K] \to [S,K] \]
inducido por
\[(g \colon S' \to K) \mapsto (g \circ f \colon S \to K)\]
es una biyección.
\end{defn}
\begin{propn} \label{whiteheadsset}
    Sean $K$ y $K'$ complejos de Kan. Un mapeo $f \colon K \to K'$ es una equivalencia débil si y solo si $f$ induce un isomorfismo en grupos de homotopía \[\pi_n(f) \colon \pi_n(K,v) \xrightarrow{\cong} \pi_n(K', f(v))\] para todo punto base $v$. 
\end{propn}
Note la analogía entre la Definición \ref{equivdebilsset} y Proposición \ref{whiteheadsset} con el Teorema \ref{Whitehead}. 

Resumamos algunas de las construcciones y enunciados que hemos discutido hasta este punto. Tenemos dos categorías $\mathsf{Top}$ y $\mathsf{sSet}$ cada una con una noción de equivalencia débil. También tenemos una clase de objetos especiales en cada categoría:
los $CW$-complejos en $\mathsf{Top}$ y los complejos de Kan en $\mathsf{sSet}$. Tenemos dos funtores entre estas categorías
\[\text{Sing} \colon \mathsf{Top} \leftrightarrows  \mathsf{sSet} \colon |- |.\] Para todo espacio topológico $X$, $\text{Sing}(X)$ es un complejo de Kan y para todo conjunto simplicial $S$, $|S|$ es un $CW$-complejo. Claramente, estos dos funtores no establecen una equivalencia (estricta) de categorías. En la siguiente discusión explicaremos el sentido en que los funtores $\text{Sing}$ y $|-|$ establecen una equivalencia de \textit{teorías de homotopía}, es decir, estos funtores resultan ser mutuamente inversos salvo equivalencias débiles naturales en cada contexto.

Para un conjunto simplicial $X$, definimos un morfismo de conjuntos simpliciales \[ \eta_X \colon X \to \Sing(|X|)\] dado en nivel $n$ por la función que envía un elemento $x\in X_n$ a la función $(f\colon |\Delta^n| \to |X|) \in \text{Sing}(|X|)_n$ definida como $f(t)=[(x,t)].$  La colección de morfismos $\eta_X$ define una transformación natural \[\eta \colon \text{id}_{\mathsf{sSet}} \Rightarrow \text{Sing} \circ | - |.\]

Dado un espacio topológico $T$, definimos un mapeo continuo \[\varepsilon_T \colon |\Sing(T)| \to T\] como la función que envía $[(f,t)]$ a $f(t)$ (donde $f\colon |\Delta^n|\to T$ y $t\in |\Delta^n|$). La colección de morfismos $\varepsilon_T$ define una transformación natural
\[ \varepsilon \colon |-| \circ \text{Sing} \Rightarrow \text{id}_{\mathsf{Top}}.\]

\begin{thm} Para todo conjunto simplcial $S$ y todo espacio topológico $X$, los morfismos naturales

\[ \eta_S \colon S \to \emph{Sing}(|S|) \text{ }\emph{   y   }\text{ } \varepsilon_X \colon |\emph{Sing}(X)| \to X\] son equivalencias débiles en $\mathsf{sSet}$ y $\mathsf{Top}$, respectivamente. En particular, todo conjunto simplicial es débil equivalente a un complejo de Kan y todo espacio topológico es débil equivalente a un CW-complejo.
\end{thm}
\subsection{Nervio topológico} Finalmente, terminamos comentando un poco mas sobre la naturaleza categórica de un complejo de Kan. En el ejemplo \ref{caminos} describimos un funtor \[\mathsf{P} \colon \mathsf{Top} \to \mathsf{Cat}_{\mathsf{Top}}\] que, para cada espacio $X$, asocia una categoría enriquecida en espacios topológicos $\mathsf{P}(X)$ tal que al tomar componentes conexos en los morfismos obtenemos el grupoide fundamental de $X$ (ejemplo \ref{grupoidefundamental}). Por otro lado, el funtor \[\text{Sing} \colon \mathsf{Top} \to \mathsf{sSet}\] asocia a cada espacio $X$ un complejo de Kan $\text{Sing}(X)$. Siguiendo la discusión previa a la Proposición \ref{SingesKan}, podemos pensar $\text{Sing}(X)$ como una categoría cuyos objetos son los puntos de $X$ (i.e. los elementos de $\text{Sing}(X)_0$), morfismos los caminos $\mathbf{\Delta}^1 \to X$ (i.e. los elementos de $\text{Sing}(X)_1$) y la regla de composición está definida salvo homotopías coherentes dadas por mapeos $\mathbf{\Delta}^n \to X$ (i.e. los elementos de $\text{Sing}(X)_n$) para $n\geq 2$. 

Una pregunta natural es cómo las construcciones $\mathsf{P}(X)$ y $\text{Sing}(X)$ están relacionadas. Se puede dar una respuesta a través del funtor \textit{nervio topológico}
\[ N_{\mathsf{Top}} \colon \mathsf{Cat}_{\mathsf{Top}} \to \mathsf{sSet} \] 
que describimos a continuación. 

Para cada entero $k \geq 0$, sea $\mathfrak{C}[k]$ un objeto en $\mathsf{sCat}$ (es decir, un funtor $\mathfrak{C}[k] \colon \Delta^{op} \to \mathsf{Cat}$) dado por
\[ \mathfrak{C}[k] \colon [n] \mapsto \mathfrak{C}[k]_n\]
donde $\mathfrak{C}[k]_n$ es la categoría con objetos $\{0,1,\ldots, k\}$ y con morfismos entre dos objetos $i,j \in \{0,1\ldots, k\}$ definidos de la siguiente manera. Si $i>j$, entonces $\mathfrak{C}[k]_n(i,j)$ es el conjunto vacío. Si $i \leq j$, definimos $\mathfrak{C}[k]_n(i,j)$ como el conjunto $N(P_{i,j})_n$ de $n$-símplices del conjunto simplicial $N(P_{i,j})$ donde $N$ denota el nervio y $P_{i,j}$ es la categoría cuyos objetos son los subconjuntos $U \subseteq \{0,1,\ldots,k\}$ tal que $i,j \in U$ y morfismos son inclusiones entre subconjuntos. La regla de composición en $\mathfrak{C}[k]_n$ está determinada por el funtor
\[P_{j,l} \times P_{i,j} \to P_{i,l}\]
inducido por la unión de subconjuntos
\[ U \times V \mapsto U \cup V.\] Dejaremos al lector la tarea de explicar como se debería definir el funtor $\mathfrak{C}[k] \colon \Delta^{op} \to \mathsf{Cat}$ a nivel de morfismos. Además note que un morfismo $\theta \colon [k] \to [m]$ en $\Delta$ induce un morfismo $\mathfrak{C}[\theta] \colon \mathfrak{C}[k] \to \mathfrak{C}[m]$ en $\mathsf{sCat}$. Podemos pensar $\mathfrak{C}[k]$ como una versión de la categoría $[k]=\{ 0 \to 1 \to \cdots \to k\}$
pero agrandada con estructura adicional que nos permitirá mantener un registro combinatorio de homotopías superiores.

Si $\mathsf{C} \in \mathsf{Cat}_{\mathsf{Top}}$, definimos $\text{Sing}(\mathsf{C}) \in \mathsf{sCat}$ declarando que cada categoría $\text{Sing}(\mathsf{C})[n]$ tenga los mismos objetos que $\mathsf{C}$ y como morfismos el conjunto de $n$-símplices $\text{Sing}(\mathsf{C}(x,y))_n$ para todo par de objetos $x$ y $y$. Los funrotores de caras
\[\text{Sing}(\mathsf{C})[n]\to \text{Sing}(\mathsf{C})[n-1]\] y degeneraciones
\[\text{Sing}(\mathsf{C})[n]\to \text{Sing}(\mathsf{C})[n+1]\] están definidos de la manera evidente: identidad en objetos y caras y degeneraciones simpliciales en los morfismos.

Definamos
$N_{\mathsf{Top}}(\mathsf{C})$ como el conjunto simplicial cuyo conjunto de $k$-símplices está dado por
\[ N_{\mathsf{Top}}(\mathsf{C})_k = \mathsf{sCat}(\mathfrak{C}[k], \text{Sing}(\mathsf{C})).\]
Para todo morfismo $\theta \colon [k] \to [m]$ obtenemos una función
\[  N_{\mathsf{Top}}(\mathsf{C})_m \xrightarrow{\mathfrak{C}[\theta]^*}  N_{\mathsf{Top}}(\mathsf{C})_k\]
así que $N_{\mathsf{Top}}(\mathsf{C})$ determina un conjunto simplicial y esta construcción define un funtor $N_{\mathsf{Top}} \colon \mathsf{Cat}_{\mathsf{Top}} \to \mathsf{sSet}$.

Denotemos por
\[ \pi_0 \colon \mathsf{Cat}_{\mathsf{Top}} \to \mathsf{Cat} \]
el funtor que aplica el funtor de componentes conexos a nivel de espacios de morfismos.

\begin{thm} 
\text{ }
\begin{enumerate}
    \item Para toda categoría enriquecida en espacios topológicos $\mathsf{C}$, $N_{\mathsf{Top}}(\mathsf{C})$ es una cuasi-categoría. Si todos los morfismos en $\pi_0(\mathsf{C})$ son isomorfismos, entonces $N_{\mathsf{Top}}(\mathsf{C})$ es un complejo de Kan. 
    \item Para todo espacio topológico $X$ existe una equivalencia débil natural de complejos de Kan
\[ \gamma_X \colon \emph{Sing}(X) \to N_{\mathsf{Top}}(\mathsf{P}(X)).\]
En particular, el espacio topológico $X$ puede ser recuperado funtorialmente, salvo equivalencia débil, de la categoría enriquecida en espacios topológicos $\mathsf{P}(X)$. 
\end{enumerate}
\end{thm}
Para más detalles sobre el nervio topológico vea \cite[\href{https://kerodon.net/tag/00KM}{Tag 00KM}]{kerodon} y el primer capítulo de \cite{htt}.
\clearpage
\pagestyle{plain}
\printbibliography
\end{document}